\documentclass[12pt]{amsart}

\usepackage[textheight=23.5cm, left=1.2cm, right=1.2cm]{geometry}

\usepackage[english]{babel}
\usepackage[T1]{fontenc}
\usepackage[utf8]{inputenc}

\theoremstyle{plain}
    \newtheorem{theorem}{Theorem}[section]
    \newtheorem{lemma}[theorem]{Lemma}
    \newtheorem{proposition}[theorem]{Proposition}

\theoremstyle{definition}

\usepackage[
    draft = false,
    unicode = true,
    colorlinks = true,
    allcolors = blue,
    hyperfootnotes = true,
    citecolor = red
]{hyperref}
\usepackage{amsfonts,amsmath,amssymb,amscd,amsthm,url}
\usepackage[inline]{enumitem}
\usepackage{graphicx,epsf,afterpage}
\usepackage{soul}
\usepackage{multicol}
\usepackage{array}
\usepackage{braket}
\usepackage{epigraph}
\usepackage{rotating}

\usepackage[
backend=biber,
style=numeric-comp,
sorting=none,
giveninits=true
]{biblatex}
\usepackage{csquotes}

\usepackage{xcolor}

\usepackage{ulem}

\usepackage{tikz}
\usetikzlibrary{positioning, calc, intersections, through, arrows, matrix, chains, math}
\usetikzlibrary{shadows}
\usetikzlibrary{graphs}
\usetikzlibrary{graphs.standard}
\usetikzlibrary{patterns}
\usetikzlibrary{decorations.pathreplacing,calligraphy,backgrounds}

\newcommand\norm[1]{\ensuremath{\left\lVert#1\right\rVert}}
\newcommand\abs[1]{\ensuremath{\left\lvert#1\right\rvert}}

\renewcommand{\Pr}{\mathrm{P}}

\DeclareMathOperator{\Expect}{\mathbb{E}}

\DeclareMathOperator{\Var}{Var}

\DeclareMathOperator{\tr}{tr}

\newcommand{\Lcal}{\mathcal{L}}
\newcommand{\Hcal}{\mathcal{H}}
\newcommand{\Tcal}{\mathcal{T}}

\newcommand{\R}{\ensuremath{\mathbb{R}}}

\newcommand{\Cx}{\ensuremath{\mathbb{C}}}

\renewcommand{\leq}{\leqslant}

\newcounter{mcnt}

\newcounter{wordcnt}

\begin{document}

\title{Limit theorems for the evolution of quantum  pure states}

\author{S.\,V. Dzhenzher and V.\,Zh. Sakbaev}

\begin{abstract}
    Limit theorems of strong law of large numbers and central limit theorem types are obtained for the compositions of independent identically distributed random unitary channels.
\end{abstract}

\thanks{\hspace{-4mm}S.\,V. Dzhenzher: sdjenjer@yandex.ru. orcid: 0009-0008-3513-4312
\\
V.\,Zh. Sakbaev: fumi2003@mail.ru. orcid: 0000-0001-8349-1738
\\
V.\,Zh. Sakbaev: Keldysh Institute of Applied Mathematics of Russian Academy of Sciences 125047, Moscow, Russia
\\
Both authors: Moscow Institute of Physics and Technology 141701, Dolgoprudny, Russia}

\maketitle
\thispagestyle{empty}

\emph{Keywords: Random unitary semigroups; Strong Law of Large Numbers; Central Limit Theorem; quantum states; Weak Operator Topology.}

\vspace{5mm}

\emph{MSC: 47D03, 60F05, 60F15, 81P16.}

\section{Introduction}

The field of random quantum channels and random quantum dynamical semigroups attracts significant interest from mathematical physics.
A general theory of random semigroups was developed in \cite{Skorokhod, Butsan-1977}.
The limit theorems for products of independent identically distributed (i.i.d.) random linear operators and semigroups of such operators are studied in \cite{Tutubalin, Furstenberg, Skorokhod}.
The limit theorems for the compositions of i.i.d. random operators have applications in differential equations \cite{Berger, Kalmetev, Zamana, Sakbaev-Smolyanov-Zamana} and in the description of the dynamics of open quantum systems \cite{Aaron, Pechen, DzhenzherSakbaev24, Dzhenzher25, AmosovSakbaev}.
Random unitary channels and random quantum dynamical semigroups naturally arise in the consideration of open quantum systems dynamics \cite{Accardi-Volovich-Lu-2002, Kempe}.
In the paper \cite{Teretenkov2023}, the approach of an open system to the Markovian dynamics is described.
The limit behavior of the sequence of quantum random walks of i.i.d. unitary channels is described in \cite{Joye, GOSS-2022, DzhenzherSakbaev25-SLLN}.

Random unitary groups arise in quantum mechanics, in particular, in the ambiguity of the quantization procedure \cite{OSS-2016}. In \cite{AmosovSakbaev}, different random semigroups in quantum states are considered. In particular, the semigroups of shift operators are discussed there. Random semigroups, and thus random quantum states, are examined from the point of view of averaging over some measure.

In this paper, we consider compositions of random unitary channels.
In particular, we obtain the Strong Law of Large Numbers (SLLN) and the Central Limit Theorem (CLT) type of convergence results for random semigroups of shifts, acting on pure non-random quantum states.
We obtain such results in terms of kernels of quantum states (theorems~\ref{t:slln-kernel} and~\ref{t:clt-kernel}), and in the Weak Operator Topology (WOT) of quantum states (theorems~\ref{t:slln-wot}, \ref{t:clt-wot}, \ref{t:random-walk}, and~\ref{t:ell1-wot}).
We also consider another semigroup of random unitary transforms (Theorem~\ref{t:analogues}).

The structure of this paper is as follows.
In \S\ref{s:basic}, we describe the main notions of quantum states, their presentation and classicication.
In \S\ref{s:coord-repr}, we obtain the SLLN and the CLT type convergence results for compositions of i.i.d. random shift operators in terms of density operator kernels.
In \S\ref{s:wot}, we define different types of convergence (almost sure, in distribution, and in $L_1$ sense) in WOT for compositions of i.i.d. random unitary channels associated with shifts on random vectors of the configuration space, and obtain sufficient conditions for above types of convergence.
In \S\ref{s:conclude}, we analyze the obtained results and discuss further research.

\section{Basic notions}\label{s:basic}

Consider the Hilbert space \(\Hcal := \Lcal_2(\R^d)\) over $\Cx$, where the integral is taken by the classical Lebesgue measure on \(\R^d\).
For \(a\in\R^d\), denote by \(S_a\colon \Hcal\to \Hcal\) the shift operator defined by
\[
    (S_a\ket{u})(x) := u(x+a).
\]
We identify $u$ and \(\ket{u}\) if necessary.

Let \((\Omega, \mathcal{F}, \mathrm{P})\) be a probability space.
It is known \cite{Orlov-Sakbaev-Shmidt-2021} that for independent identically distributed (i.i.d.) random vectors \(\xi_1, \ldots, \xi_n\) in $\R^d$ with finite expected value and variance the compositions of independent shifts \(S_{\sqrt{\frac{t}{n}}\xi_1}\ldots S_{\sqrt{\frac{t}{n}}\xi_n}\ket{u}\) converge in distribution to a convolution of $\ket{u}$ and the Gaussian measure.
In this text we explore the composition of shifts acting on a quantum state \(\ket{u}\bra{u}\).

In this text we consider \textbf{\emph{normal quantum states}}, that are nucleus operators from the space $\Tcal_1(\Hcal)$; in other words, they are operators with trace.
Note that the dual \(\Tcal_1(\Hcal)^* = B(\Hcal)\) is the space of all bounded linear operators on $\Hcal$.
The second dual space \(\Tcal_1(\Hcal)^{**} = B(\Hcal)^*\) is the space of all quantum states; but this general case will not arise in this paper.
So we will shortly call the normal quantum states just as <<quantum states>>.

For any $\ket{u} \in \Hcal$ we denote by \(\rho [u] := \ket{u}\bra{u} \in \Tcal_1(\Hcal)\) the \textbf{\emph{pure state}}.

\section{Coordinate representation}\label{s:coord-repr}

For any quantum state $\rho\in\Tcal_1(\Hcal)$ and any $x,y \in \R^d$ denote by
\(\rho(x,y) \in \R\)
its \emph{kernel}; that is, the function for any $u,v\in\Hcal$ given by
\[
    \bra{u}\rho\ket{v} = \int \rho(x,y)u(x)\overline{v(y)}\,dxdy.
\]
In Fourier images, we will write $\alpha,\beta$ instead of $x,y$.

Denote by $F\colon\Hcal\to\Hcal$ the Fourier transform.
We identify \(F\ket{u}\) and \(\ket{Fu}\).

Denote by $\Expect$ the expected value of a random variable.

The following two theorems generalize the SLLN and the CLT for kernels of quantum states.

\begin{theorem}[SLLN for kernel]\label{t:slln-kernel}
    Let \(\ket{u}\in \Hcal\).
    Let \(\{\xi_n\}\) be i.i.d. random vectors with finite expectation $\mu\in\R^d$.
    Then for any \(x,y \in \R^d\)
    \[
        \rho\left[S_{\xi_1/n}\ldots S_{\xi_n/n} u\right](x, y)
        \xrightarrow[n\to\infty]{a.s.}
        \rho[S_{\mu}u](x, y).
    \]
\end{theorem}

\begin{theorem}[CLT for kernel]\label{t:clt-kernel}
    Let \(\ket{u}\in \Hcal\).
    Let \(\{\xi_n\}\) be i.i.d. zero-mean random vectors with the covariance matrix~$\Sigma$.
    Then for any \(x,y \in \R^d\)
    \[
        \rho \left[S_{\xi_1/\sqrt{n}}\ldots S_{\xi_n/\sqrt{n}} u\right](x, y)
        \xrightarrow[n\to\infty]{d}
        \rho[S_{\eta}u](x, y),
    \]
    where \(\eta\sim N(0,\Sigma)\) is a random Gaussian vector.
    In particular,
    \[
        \Expect F\rho \left[S_{\xi_1/\sqrt{n}}\ldots S_{\xi_n/\sqrt{n}} u\right] F^{-1}(\alpha, \beta)
        \to
        e^{-\frac{1}{2} (\alpha-\beta)^T\Sigma(\alpha-\beta)}\rho[Fu](\alpha, \beta)
        \quad\text{as}\quad n\to\infty.
    \]
\end{theorem}

These theorems follow from Continuous Mapping Theorem~\ref{t:cmt} and the following proposition.

\begin{proposition}
\label{p:ft-shift}
    Let \(\ket{u}\in \Hcal\) and $a \in \R^d$.
    Then for any $\alpha,\beta \in \R^d$
    \[
        F\rho\left[S_au\right]F^{-1}(\alpha,\beta) = e^{i(a, \alpha-\beta)}\rho\left[Fu\right](\alpha, \beta).
    \]
\end{proposition}
\begin{proof}
    This follows since for any \(\ket{v} \in \Hcal\) we have
    \[
        F\rho[v]F^{-1} = F\ket{v}\bra{v}F^{-1} = \ket{Fv}\bra{Fv} = \rho[Fv],
    \]
    and since
    \[
        \left(FS_{a}\ket{u}\right)(\alpha) = e^{i(a, \alpha)}F\ket{u}(\alpha).
    \]
\end{proof}

\begin{theorem}[Continuous Mapping {\cite{Billingsley2013convergence}}]\label{t:cmt}
    Let \(\mathcal{S}\) be a metric space.
    Let \(\{\xi_n\}\) be random variables on \(\mathcal{S}\), converging to \(\xi\) either almost surely, or in probability, or in distribution.
    Let \(f\colon\mathcal{S}\to\R\) be a Borel function, which is continuous \(\Pr_\xi\)-almost surely.
    Then random variables \(f(\xi_n)\) converge to a random variable \(f(\xi)\) in the same sense.
\end{theorem}

\begin{proof}[Proof of Theorem~\ref{t:slln-kernel}]
    Since shifts form a one-parameter group, by Proposition~\ref{p:ft-shift} we have
    \[
        F\rho\left[S_{\xi_1/n}\ldots S_{\xi_n/n} u\right] F^{-1}(\alpha, \beta)
        = e^{i\left(\frac{\sum \xi_i}{n}, \alpha-\beta\right)}\rho[Fu](\alpha, \beta).
    \]
    In the right part, only the complex exponent is random.
    By the Strong Law of Large Numbers for random vectors,
    \[
        \frac{\sum\xi_i}{n} \xrightarrow{a.s.} \mu.
    \]
    So by Continuous Mapping Theorem~\ref{t:cmt} we obtain that
    \[
        F\rho\left[S_{\xi_1/n}\ldots S_{\xi_n/n} u\right] F^{-1}(\alpha, \beta)
        \xrightarrow{a.s.}
        e^{i\left(\mu, \alpha-\beta\right)}\rho[Fu](\alpha, \beta).
    \]
    Now the result follows by acting with the inverse Fourier transform in Proposition~\ref{p:ft-shift}.
\end{proof}

\begin{proof}[Proof of Theorem~\ref{t:clt-kernel}]
    Since shifts form a one-parameter group, by Proposition~\ref{p:ft-shift} we have
    \[
        F\rho\left[S_{\xi_1/\sqrt{n}}\ldots S_{\xi_n/\sqrt{n}} u\right] F^{-1}(\alpha, \beta)
        = e^{i\left(\frac{\sum \xi_i}{\sqrt{n}}, \alpha-\beta\right)}\rho[Fu](\alpha, \beta).
    \]
    In the right part only the complex exponent is random.
    By the Central Limit Theorem for random vectors,
    \[
        \frac{\sum\xi_i}{\sqrt{n}} \xrightarrow{d} \eta.
    \]
    So by Continuous Mapping Theorem~\ref{t:cmt} we obtain that
    \[
        F\rho \left[S_{\xi_1/\sqrt{n}}\ldots S_{\xi_n/\sqrt{n}} u\right] F^{-1}(\alpha, \beta)
        \xrightarrow{d}
        e^{i\left(\eta, \alpha-\beta\right)}\rho[Fu](\alpha, \beta).
    \]
    Now the result follows by acting with the inverse Fourier transform in Proposition~\ref{p:ft-shift}.
    The particular case is obtained by definition of convergence in distribution, since
    \[
        \Expect e^{i(a,\eta)} = e^{-\frac{1}{2}a^T\Sigma a}
    \]
    is the characteristic function of the Gaussian vector.
\end{proof}

\section{Weak Operator Topology of quantum states}\label{s:wot}

Here we develop another approach for convergence in WOT of quantum states.
By \cite[Theorem~2.7.2]{Holevo2011-ch2} this topology of linear continuous functionals can be described by bounded operators \(A\in B(\Hcal)\) and functionals
\[
    \rho \mapsto \tr(\rho A).
\]
Note that apart from the choice of the topology we need to choose the probabilistic mode of convergence.
We will work with three of them: almost sure, in distribution, and in \(L_1(\Pr)\).
So, when we say <<random states \(\rho_n\) converge in some probabilistic mode in WOT to \(\rho\)>>, we will mean that for any \(A\in B(\Hcal)\) the random variables \(\tr (\rho_n A)\) converge in that probabilistic mode to the random variable \(\tr(\rho A)\).

\begin{theorem}[SLLN in WOT]
\label{t:slln-wot}
    Let \(\ket{u}\in \Hcal\).
    Let \(\{\xi_n\}\) be i.i.d. random vectors with finite expectation $\mu\in\R^d$.
    Then the random pure states \(\rho\left[S_{\xi_1/n}\ldots S_{\xi_n/n} u\right]\) converge a.s. in WOT to the pure state \(\rho\left[S_{\mu}u\right]\);
    in other words,
    for any \(A \in B(\Hcal)\)
    \[
        \tr\left(\rho\left[S_{\xi_1/n}\ldots S_{\xi_n/n} u\right]A\right) \xrightarrow[n\to\infty]{a.s.} \tr\Bigl(\rho\left[S_{\mu}u\right]A\Bigr).
    \]
\end{theorem}

\begin{theorem}[CLT in WOT]
\label{t:clt-wot}
    Let \(\ket{u}\in \Hcal\).
    Let \(\{\xi_n\}\) be i.i.d. zero-mean random vectors with the covariance matrix $\Sigma\in\R^{d\times d}$.
    Then the random pure states \(\rho\left[S_{\xi_1/\sqrt{n}}\ldots S_{\xi_n/\sqrt{n}} u\right]\) converge in distribution in WOT to the random pure state \(\rho\left[S_{\eta}u\right]\), where \(\eta \sim N(0, \Sigma)\) is a random Gaussian vector;
    in other words,
    for any \(A \in B(\Hcal)\)
    \[
        \tr\left(\rho\left[S_{\xi_1/\sqrt{n}}\ldots S_{\xi_n/\sqrt{n}} u\right]A\right) \xrightarrow[n\to\infty]{d} \tr\Bigl(\rho\left[S_{\eta}u\right]A\Bigr).
    \]
    Moreover, we have
    \[
        \Expect\tr\left(\rho\left[S_{\xi_1/\sqrt{n}}\ldots S_{\xi_n/\sqrt{n}} u\right]A\right) \xrightarrow[n\to\infty]{} \Expect\tr\Bigl(\rho\left[S_{\eta}u\right]A\Bigr).
    \]
\end{theorem}

\begin{theorem}[Random walk]\label{t:random-walk}
    Let \(d=1\).
    Let \(\xi_{n,k}\) be independent random variables with \(\Expect\xi_{n,k}=0\), \(\Var \xi_{n,k} = \frac{1}{n}\), and
    \[
        \max_{1\leq k\leq n} \Pr\{\abs{\xi_{n,k}}>\varepsilon\} \xrightarrow[n\to\infty]{} 0.
    \]
    Take any \(u \in \Hcal\) and consider the <<random polygonal line>> \(\zeta_n(t)\) constructed by the points 
    \[
        \left(\frac{k}{n}, \rho\left[S_{\xi_{n,1}}\ldots S_{\xi_{n,k}}u\right]\right)
    \]
    for \(k=1,\ldots,n\).
    Let \(w(t)\) be the standart Wiener process.
    
    Then \(\xi_n\) converges in distribution in WOT to the random states-valued process \(\rho[S_{w}u]\);
    in other words, for any \(A \in B(\Hcal)\) the random process
    \[
        \tr\left(\xi_n(t)A\right)
    \]
    converges in distribution to
    \[
        \tr\left(\rho[S_{w(t)}u]A\right)
    \]
\end{theorem}

All the last three theorems follow from Continuous Mapping Theorem~\ref{t:cmt} and from Lemma~\ref{l:tr-cont}.
The <<moreover>> part in Theorem~\ref{t:clt-wot} follows since random shifts of the pure state \(\ket{u}\) still lie in the ball of the radius \(\norm{u}\), and hence all random variables there are equi-bounded.

\begin{lemma}[{\cite[Lemma~3.1]{AmosovSakbaev}}]\label{l:trace-norm}
    Let $\Hcal$ be a Hilbert space, and \(u,v \in \Hcal\) be normed.
    Then
    \[
        \tr \abs{\rho[u] - \rho[v]} \leq 2\norm{u-v}.
    \]
\end{lemma}

\begin{lemma}\label{l:tr-cont}
    For any \(u \in \Hcal\) and \(A \in B(\Hcal)\)
    the function \( x \mapsto \tr\left(\rho\left[S_xu\right]A\right)\) is a continuous function \( \R^d \to \R\).
\end{lemma}

\begin{proof}
    This holds since by Lemma~\ref{l:trace-norm} for some constant \(C>0\)
    \[
        \abs{\tr\left(\rho\left[S_x u\right]A\right) -\tr\left(\rho\left[S_y u\right]A\right)}
        \leqslant
        \tr\abs{\rho\left[S_x u\right] -\rho\left[S_y u\right]} \cdot \norm{A} \leqslant
        C \norm{S_x u - S_y u}
    \]
    and since (by the Plancherel theorem)
    \[
        \norm{S_xu - S_yu} = \norm{FS_xu - FS_yu} = \abs{e^{ix} - e^{iy}} \norm{u} \xrightarrow[x \to y]{}0.
    \]
\end{proof}

\begin{theorem}\label{t:ell1-wot}
    Let \(\ket{u}\in \Hcal\).
    Let \(\{\xi_n\}\) be i.i.d. random vectors with expectation $\mu$.
    Then the random pure states \(\rho\left[S_{\xi_1/n}\ldots S_{\xi_n/n} u\right]\) converge in WOT in $L_1(\Pr)$ to the pure state \(\rho\left[S_\mu u\right]\); in other words, 
    for any \(A \in B(\Hcal)\)
    \[
        \Expect \abs{\tr\left(\rho\left[S_{\xi_1/n}\ldots S_{\xi_n/n} u\right]A\right) - \tr\left(\rho\left[S_\mu u\right]A\right)}
        \xrightarrow[n\to\infty]{} 0.
    \]
    %
\end{theorem}

\begin{proof}
    Since the shift does not change the norm, all random states \(S_{\xi_1/n}\ldots S_{\xi_n/n} u\) lie in the ball of the radius \(\norm{u}\).
    Hence the theorem follows by Theorem~\ref{t:slln-wot} and by \cite[Lebesgue's Dominated Convergence Theorem~16.4]{Billingsley2012probability}.
    %
\end{proof}

Let us now consider for any \(a \in \R^d\) the <<impulse>> operator \(R_a \colon \Hcal\to\Hcal\) defined by
\[
    (R_a \ket{u})(x) := e^{iax}u(x).
\]

\begin{theorem}\label{t:analogues}
    The analogues of theorems~\ref{t:slln-wot}, \ref{t:clt-wot}, \ref{t:random-walk} and~\ref{t:ell1-wot} hold for \(R_a\) instead of \(S_a\).
\end{theorem}

This again holds by Continuous Mapping Theorem~\ref{t:cmt} and the following lemma.

\begin{lemma}
    For any \(u \in \Hcal\) and \(A \in B(\Hcal)\)
    the function \( x \mapsto \tr\left(\rho\left[R_xu\right]A\right)\) is a continuous function \( \R^d \to \R\).
\end{lemma}

\begin{proof}
    It is clear that \(FS_a = R_aF\), and so \(R_a = FS_aF^{-1}\).
    Now the lemma follows again by Lemma~\ref{l:trace-norm} since for some $C>0$ depending on $u,A$ we have
    \begin{multline*}
        \abs{\tr\left(\rho\left[R_x u\right]A\right) -\tr\left(\rho\left[R_y u\right]A\right)}
        \leqslant \\ \leq
        \tr\abs{\rho\left[R_x u\right] -\rho\left[R_y u\right]} \cdot \norm{A} \leqslant
        C \norm{R_x u - R_y u}  = C\abs{e^{ix} - e^{iy}} \norm{F^{-1}u} \xrightarrow[x \to y]{}0.
    \end{multline*}
\end{proof}

\section{Conclusions}\label{s:conclude}

The analogues of theorems~\ref{t:slln-wot}, \ref{t:clt-wot}, \ref{t:random-walk}, and~\ref{t:ell1-wot} may not hold for some other measure instead of the classical Lebesgue measure on~$\R^d$.
Indeed, in the proof of Lemma~\ref{l:tr-cont}, it is used that the Fourier transform is an isometry (in other words, the Plancherel theorem is used);
but, it does not hold for the general measure~$m$.
Nevertheless, the analogues of these theorems hold for continuous functions \(u\in\Hcal\) and for finite $\sigma$-additive measures, that can be proved by the Egoroff's theorem.
Generally speaking, these analogues may hold for some other measure, but the proof would require a technique which is different from the one described here, since we use greatly the properties of the Fourier transform (in particular, the Plancherel equality).
It would be interesting to obtain the analogues of the mentioned theorems for the other measures.

It is also interesting to obtain the analogues of these theorems for the Hilbert space \(\Hcal\) of functions on some infinite-dimensional space instead of $\R^d$. The case of \(\ell_2(\R)\) with some analogue of the Lebesgue measure looks like the most obvious choice to consider.

Finally, we have discussed here only normal quantum states.
The case of singular quantum states seems to be interesting case to develop.
In particular, it is interesting to consider some other random transforms instead of random shifts, so that random states would not be pure.
For example, it can be random unitary transforms with some special random generators.

\printbibliography

\end{document}